\documentclass[fleqn,preprint,3p,a4paper]{elsarticle}

\usepackage{amssymb}
\usepackage{amsmath}
\usepackage{amsthm}
\usepackage{dcolumn}
\usepackage{endnotes}
\usepackage{tabularx}
\usepackage[matrix,arrow]{xy}
\usepackage{wasysym}

\newtheorem{theorem}{Theorem}[section]
\newtheorem{proposition}[theorem]{Proposition}
\newtheorem{lemma}[theorem]{Lemma}
\newtheorem{corollary}[theorem]{Corollary}

\theoremstyle{definition}
\newtheorem{definition}[theorem]{Definition}
\newtheorem{example}[theorem]{Example}
\newtheorem{remark}[theorem]{Remark}

\newtheorem{problem}[theorem]{Problem}

\newcommand{\II}{\begin{enumerate}}
\newcommand{\III}{\end{enumerate}}

\newcommand{\ir}{{\mathsf{Irr}}}

\newcommand{\cl}{{\rm cl}}
\newcommand{\ii}{{\rm int}}

\newcommand{\ua}{\mathord{\uparrow}}
\newcommand{\da}{\mathord{\downarrow}}

\newcommand{\mk}{\mathord{\mathsf{K}}}
\newcommand{\wdd}{\mathord{\mathsf{WD}}}

\newcommand{\kf}{\mathord{\mathsf{RD}}}

\journal{Topology and its applications}

\begin{document}

\begin{frontmatter}



\title{On topological Rudin's lemma, well-filtered spaces and sober spaces\tnoteref{t1}}
\tnotetext[t1]{This research was supported by the National Natural Science Foundation of China (No. 11661057); the Natural Science Foundation of Jiangxi Province , China (No. 20192ACBL20045); and NIE ACRF (RI 3/16 ZDS), Singapore}

\author[X. Xu]{Xiaoquan Xu\corref{mycorrespondingauthor}}
\cortext[mycorrespondingauthor]{Corresponding author}
\ead{xiqxu2002@163.com}
\address[X. Xu]{School of Mathematics and Statistics,
Minnan Normal University, Zhangzhou 363000, China}
\author[D. Zhao]{Dongsheng Zhao}
\address[D. Zhao]{Mathematics and Mathematics Education,
National Institute of Education Singapore, \\
Nanyang Technological University,
1 Nanyang Walk, Singapore 637616}
\ead{dongsheng.zhao@nie.edu.sg}

\begin{abstract}
Based on topological Rudin's Lemma, we investigate two new kinds of sets - Rudin sets and well-filtered determined sets in $T_0$ topological spaces. Using such sets, we formulate and prove some new characterizations for well-filtered spaces and sober spaces. Part of the work was inspired by Xi and Lawson's work on well-filtered spaces.
Our study also lead to a new class of spaces - strong $d$-spaces and some problems whose solutions will strengthen our understanding of the related structures.
\end{abstract}





\end{frontmatter}



In domain theory, the $d$-spaces, well-filtered spaces and sober spaces form three of the most important classes of spaces. Rudin's Lemma has played a crucial role in studying such spaces. The original application of Rudin's Lemma was in answering some questions on quasicontinuous dcpos. In recent years, it has been used to study the various aspects of well-filtered spaces, initiated by Heckmann and Keimel
\cite{Klause-Heckmann}. In this paper, inspired by the topological version of Rudin's Lemma by Heckmann and Keimel, Xi and Lawson's work \cite{Xi-Lawson-2017} on well-filtered spaces and our recent work \cite{Shenchon} on the well-filtered reflections of $T_0$ spaces, we investigate two new kinds of sets - Rudin sets and well-filtered determined sets in $T_0$ topological spaces, and use them to establish a series new characterizations of well-filtered spaces and sober spaces. Our study also leads to a new class of spaces - strong $d$-spaces, and a number of problems, whose answering will deepen our understanding of the related spaces and structures.

\section{Preliminary}

For a poset $P$ and $A\subseteq P$, let
$\mathord{\downarrow}A=\{x\in P: x\leq  a \mbox{ for some }
a\in A\}$ and $\mathord{\uparrow}A=\{x\in P: x\geq  a \mbox{
	for some } a\in A\}$. For  $x\in P$, we write
$\mathord{\downarrow}x$ for $\mathord{\downarrow}\{x\}$ and
$\mathord{\uparrow}x$ for $\mathord{\uparrow}\{x\}$.  A subset $A$
is called a \emph{lower set} (resp., an \emph{upper set}) if
$A=\mathord{\downarrow}A$ (resp., $A=\mathord{\uparrow}A$). Let $P^{(<\omega)}=\{F\subseteq P : F \mbox{~is a nonempty finite set}\}$, $P^{(\leqslant\omega)}=\{F\subseteq P : F \mbox{~is a nonempty countabel set}\}$ and $\mathbf{Fin} ~P=\{\uparrow F : F\in P^{(<\omega)}\}$. The set of all nonempty upper subsets of $P$ is denoted by $\mathbf{up}(P)$. For a nonempty subset $B$ of $P$, let $max (B)=\{b\in B : b \mbox{~ is a maximal element of~} B\}$ and $min (B)=\{b\in B : b \mbox{~ is a minimal element of~} B\}$. A
nonempty subset $D$ of $P$ is \emph{directed} if every two
elements in $D$ have an upper bound in $D$. The set of all directed sets of $P$ is denoted by $\mathcal D(P)$. A subset $I\subseteq P$ is called an \emph{ideal} of $P$ if $I$ is a directed  and a lower set. Let $\mathrm{Id} (P)$ be the poset (with the order of set inclusion) of all ideals of $P$. Dually, we define the concept of \emph{filters} and denote the poset of all filters of $P$ by $\mathrm{Filt}(P)$. $P$ is called a
\emph{directed complete poset}, or \emph{dcpo} for short, if for any
$D\in \mathcal D(P)$, $\bigvee D$ exists in $P$.

As in \cite{redbook}, the \emph{lower topology} on a poset $P$, generated
by the complements of the principal filters of $P$, is denoted by $\omega (P)$. Dually, define the \emph{upper topology} on $P$ and denote it by $\upsilon (P)$. A subset $U$ of $P$ is \emph{Scott open} if
(i) $U=\mathord{\uparrow}U$ and (ii) for any directed subset $D$ for
which $\bigvee D$ exists, $\bigvee D\in U$ implies $D\cap
U\neq\emptyset$. All Scott open subsets of $P$ form a topology,
and we call this topology  the \emph{Scott topology} on $P$ and
denote it by $\sigma(P)$. The space $\Sigma~\!\! P=(P,\sigma(P))$ is called the
\emph{Scott space} of $P$. The common refinement $\sigma(P)$ and $\omega (P)$ is called the \emph{Lawson topology}
and is denoted by $\lambda (P)$. The upper sets
form the (\emph{upper}) \emph{Alexandroff topology} $\alpha (P)$.

A poset $Q$ is called an (inf) \emph{semilattice} if for any two elements
$a, b\in Q$, $inf\{a, b\}=a\wedge b$ exists in $Q$. Dually, $Q$ is a \emph{sup semilattice} if for any two elements $a, b\in Q$, $sup\{a, b\}=a\vee b$ exists in $Q$.

\begin{definition}\label{bounded copmplete} Let $P$ be a poset.
\begin{enumerate}[\rm (i)]
\item $P$ is called a \emph{complete semilattice} if $P$ is a dcpo and every nonempty subset $P$
has an inf. In particular, a complete semilattice has a smallest
element, the inf of $P$.
\item  $P$  is called \emph{bounded complete} if every subset that is bounded above
has a sup (i.e., the least upper bound). In particular, a bounded complete poset has a smallest
element, the least upper bound of the empty set.
\end{enumerate}
\end{definition}

It is easy to see that $P$ is a complete semilattice iff $P$ is a bounded complete dcpo (see, e,g., \cite[Proposition O-2.2]{redbook}).

For a $T_0$ space $X$, we use $\leq_X$ to represent the \emph{specialization order} of $X$, that is, $x\leq_X y$ if{}f $x\in \overline{\{y\}}$. In the following, when a $T_0$ space $X$ is considered as a poset, the order always refers to the specialization order if no other explanation. Let $\mathcal O(X)$ (resp., $\mathcal C(X)$) be the set of all open subsets (resp., closed subsets) of $X$, and let $\mathcal S^u(X)=\{\ua x : x\in X\}$, $\mathcal S_c(X)=\{\overline{{\{x\}}} : x\in X\}$ and $\mathcal D_c(X)=\{\overline{D} : D\in \mathcal D(X)\}$. A space $X$ is called a \emph{$d$-space} (or \emph{monotone convergence space}) if $X$ (endowed with the specialization order) is a dcpo
 and $\mathcal O(X) \subseteq \sigma(X)$ (cf. \cite{redbook}\cite{Wyler}).

\begin{proposition}\label{d-spaceCharact} For a $T_0$ space $X$, the following conditions are equivalent:
\begin{enumerate}
            \item $X$ is a $d$-space.
            \item $\mathcal D_c(X)=\mathcal S_c(X)$.
            \item  For any $D\in \mathcal D(X)$ and $U\in \mathcal O(X)$, $\bigcap_{d\in D}\ua d\subseteq U$ implies $\ua d \subseteq U$ (i.e., $d\in U$) for some $d\in D$.
            \item  For any $D\in \mathcal D(X)$ and $A\in \mathcal C(X)$, if $D\subseteq A$, then $A\cap\bigcap_{d\in D}\ua d\neq\emptyset$.
            \item  For any $D\in \mathcal D(X)$, $\overline{D}\cap\bigcap_{d\in D}\ua d\neq\emptyset$.
\end{enumerate}
\end{proposition}
\begin{proof} (1) $\Leftrightarrow$ (2): If $X$ is a $d$-space, then for any $D\in \mathcal{D}(X)$, $\overline{D}=\overline\{sup\,D\}$, thus (1) $\Rightarrow$ (2). Conversely, if condition (2) holds, then for each $D\in \mathcal D(X)$ and $A\in \mathcal C(X)$ with $D\subseteq A$, there is $x\in X$ such that $\overline{D}=\overline{\{x\}}$, and consequently, $\bigvee D=x$ and $\bigvee D \in A$ since $\overline{D}\subseteq A$. Thus $X$ is a dcpo and $\mathcal O(X)\subseteq \sigma(X)$, hence $X$ is a $d$-space.

(1) $\Rightarrow$ (3): Since $X$ is a $d$-space, $\ua \bigvee D=\bigcap_{d\in D}\ua d\subseteq U\in \sigma(X)$. Therefore, $\bigvee D\in U$, thus
$d\in U$ for some $d\in D$.

(3) $\Rightarrow$ (4): If $A\cap\bigcap_{d\in D}\ua d=\emptyset$, then $\bigcap_{d\in D}\ua d \subseteq X\setminus A$. By condition (3), $\ua d\subseteq X\setminus A$ for some $d\in D$, which contradicts $D\subseteq A$.

(4) $\Rightarrow$ (5): Trivial.

(5) $\Rightarrow$ (1): For each $D\in \mathcal D(X)$ and $A\in \mathcal C(X)$ with $D\subseteq A$, by condition (5), $\overline{D}\cap\bigcap_{d\in D}\ua d\neq\emptyset$. Choose one $x\in \overline{D}\cap\bigcap_{d\in D}\ua d$. Then $D\subseteq \da x\subseteq \overline{D}$, hence $\overline{D}=\da x$ and $\bigvee D=x$. Therefore, $\bigvee D\in A$ because $\overline{\{\bigvee D\}}=\overline{D}\subseteq A$. Thus $X$ is a $d$-space.
\end{proof}

\begin{lemma}\label{d-spacemax}
If $X$ is a $d$-space and $A$ is a nonempty closed subset of $X$, then  $max (A)\neq\emptyset$.
\end{lemma}
\begin{proof}
By Zorn's Lemma there is a maximal chain $C$ in $A$. Since $X$ is a $d$-space, $c=\bigvee C$ exists and $c\in A$. By the maximality of $C$, we have $c\in max (A)$.
\end{proof}

A nonempty subset $A$ of a $T_0$ space $X$ is said to be \emph{irreducible} if for any $\{F_1, F_2\}\subseteq \mathcal C(X)$, $A \subseteq F_1\cup F_2$ implies $A \subseteq F_1$ or $A \subseteq  F_2$.  Denote by $\ir(X)$ (resp., $\ir_c(X)$) the set of all irreducible (resp., irreducible closed) subsets of $X$. Every subset of $X$ that is directed under $\leq_X$ is irreducible. A space $X$ is called \emph{sober}, if for any  $F\in\ir_c(X)$, there is a unique point $a\in X$ such that $F=\overline{\{a\}}$.

The following two lemmas on irreducible sets are well-known that will be used in the sequel.

\begin{lemma}\label{irrsubspace}
Let $X$ be a space and $Y$ a subspace of $X$. Then the following conditions are equivalent for a
subset $A\subseteq Y$:
\begin{enumerate}[\rm (1)]
	\item $A$ is an irreducible subset of $Y$.
	\item $A$ is an irreducible subset of $X$.
	\item ${\rm cl}_X A$ is an irreducible subset of $X$.
\end{enumerate}
\end{lemma}

\begin{lemma}\label{irrimage}
	If $f : X \longrightarrow Y$ is continuous and $A\in\ir (X)$, then $f(A)\in \ir (Y)$.
\end{lemma}

\begin{remark}\label{subspaceirr}  If $Y$ is a subspace of a space $X$ and $A\subseteq Y$, then by Lemma \ref{irrsubspace}, $\ir (Y)=\{B\in \ir(X) : B\subseteq Y\}\subseteq \ir (X)$ and  $\ir_c (Y)=\{B\in \ir(X) : B\in \mathcal C(Y)\}\subseteq \ir (X)$. If $Y\in \mathcal C(X)$, then $\ir_c(Y)=\{C\in \ir_c(X) : C\subseteq Y\}\subseteq \ir_c (X)$.
\end{remark}

For any topological space $X$, $\mathcal G\subseteq 2^{X}$ and $A\subseteq X$, let $\Diamond_{\mathcal G} A=\{G\in \mathcal G : G\bigcap A\neq\emptyset\}$ and $\Box_{\mathcal G} A=\{G\in \mathcal G : G\subseteq  A\}$. The sets $\Diamond_{\mathcal G} A$ and $\Box_{\mathcal G} A$ will be simply written as $\Diamond A$ and $\Box A$, respectively, if there is no confusion. The \emph{upper Vietoris topology} on $\mathcal{G}$ is the topology that has $\{\Box_{\mathcal{G}} U : U\in \mathcal O(X)\}$ as a base and the resulting space is denoted by $P_S(\mathcal{G})$. The \emph{lower Vietoris topology} on $\mathcal{G}$ is the topology that has $\{\Diamond U : U\in \mathcal O(X)\}$ as a subbase and the resulting space is denoted by $P_H(\mathcal{G})$. If $\mathcal{G}\subseteq \ir (X)$, then $\{\Diamond_{\mathcal{G}} U : U\in \mathcal O(X)\}$ is a topology on $\mathcal{G}$. The space $P_H(\mathcal{C}(X)\setminus \{\emptyset\})$ is called the \emph{Hoare power space} or \emph{lower space} of $X$ and is denoted by $P_H(X)$ for short (cf. \cite{Schalk}). Clearly, $P_H(X)=(\mathcal{C}(X)\setminus \{\emptyset\}, \upsilon(\mathcal{C}(X)\setminus \{\emptyset\}))$. So $P_H(X)$ is always sober (see, e.g., \cite[Corollary 4.10]{ZhaoHo}). The space $P_H(\ir_c(X))$, shortly denoted by $X^s$, with the topological embedding $\eta_X (=x\mapsto \overline{\{x\}}): X \longrightarrow P_H(\ir_c(X))$, is the \emph{canonical soberification} of $X$ (cf. \cite{redbook}).

A subset $A$ of a space $X$ is called \emph{saturated} if $A$ equals the intersection of all open sets containing it (equivalently, $A$ is an upper set with respect to the specialization order). We shall use $\mathord{\mathsf{K}}(X)$ to
denote the set of all nonempty compact saturated subsets of $X$. $X$ is called \emph{coherent} if the intersection of any two compact saturated sets is again compact. $X$ is said to be \emph{well-filtered} if it is $T_0$, and for any open set $U$ and filtered family $\mathcal{K}\subseteq \mathord{\mathsf{K}}(X)$, $\bigcap\mathcal{K}{\subseteq} U$ implies $K{\subseteq} U$ for some $K{\in}\mathcal{K}$. The space $P_S(\mathord{\mathsf{K}}(X))$, denoted shortly by $P_S(X)$, is called the \emph{Smyth power space} or \emph{upper space} of $X$ (cf. \cite{Heckmann}\cite{Schalk}). The space $P_S(\mathbf{up}(X))$ is called the \emph{Alexandroff power space}. It is easy to see that $P_S(X)$ is a subspace of $P_S(\mathbf{up}(X))$, and the specialization orders on $P_S(\mathbf{up}(X))$ is the \emph{Smyth preorder}, that is, for $K_1,K_2\in \mathbf{up}(X)$, $K_1\leq_{P_S(\mathbf{up}(X))}K_2$ if{}f $K_2\subseteq K_1$. The \emph{canonical mapping} $\xi_X: X\longrightarrow P_S(X)$, $x\mapsto\ua x$, is an order and topological embedding (cf. \cite{Heckmann}\cite{Klause-Heckmann}\cite{Schalk}). Clearly, $P_S(\mathcal S^u(X))$ is a subspace of $P_S(X)$ and $X$ is homeomorphic to $P_S(\mathcal S^u(X))$.

By Lemma \ref{irrsubspace} and Lemma \ref{irrimage}, we have the following corollary.

 \begin{corollary}\label{X-Smyth-irr} Let $X$ be a $T_0$ space and $\mathcal{A}\subseteq \mk(X)$. Then the following three conditions are equivalent:
 \begin{enumerate}[\rm (1)]
	\item $\mathcal{A}\in\ir (P_S(X))$.
	\item $\mathcal{A}\in\ir (P_S(\mathbf{up}(X)))$.
    \item $\cl_{P_S(\mathbf{up}(X))}\mathcal{A}\in \ir_c (P_S(\mathbf{up}(X)))$.
\end{enumerate}
\end{corollary}

\begin{remark}\label{meet-in-Smyth} Let $X$ be a $T_0$ space and $\mathcal A\subseteq \mk (X)$ (resp., $\mathcal A\subseteq \mathbf{up} (X)$). Then $\bigcap \mathcal A=\bigcap \overline{\mathcal A}$, here the closure of $\mathcal A$ is taken in $P_S(X)$ (resp., in $P_S(\mathbf{up}(X))$). Clearly, $\bigcap \overline{\mathcal A}\subseteq\bigcap \mathcal A$. On the other hand, for $K\in \overline{\mathcal A}$ and $U\in \mathcal O(X)$ with $K\subseteq U$ (that is, $K\in \Box U$), we have $\mathcal A\bigcap\Box U\neq\emptyset$, and hence there is a $K_U\in \mathcal A\bigcap\Box U$. Therefore $K=\bigcap \{U\in \mathcal O(X) : K\subseteq U\}\supseteq\bigcap \{K_U : U\in \mathcal O(X) \mbox{ and } K\subseteq U\}\supseteq\bigcap \mathcal A$. It follows that $\bigcap \overline{\mathcal A}\supseteq\bigcap \mathcal A$. Thus $\bigcap \mathcal A=\bigcap \overline{\mathcal A}$.
\end{remark}

\section{Rudin sets and well-filtered determined sets}

Rudin's Lemma plays a crucial role in domain theory (see [1-7]). Rudin \cite{Rudin} proved her lemma by transfinite methods. Later, Heckmann and Keimel \cite{Klause-Heckmann} established  the following topological variant of Rudin's Lemma.

\begin{lemma}\label{t Rudin} \emph{(Topological Rudin's Lemma)} Let $X$ be a topological space and $\mathcal{A}$ an
irreducible subset of the Smyth power space $P_S(X)$. Then every closed set $C {\subseteq} X$  that
meets all members of $\mathcal{A}$ contains an minimal irreducible closed subset $A$ that meets all
members of $\mathcal{A}$.
\end{lemma}

Applying Lemma \ref{t Rudin} to the Alexandroff topology on a poset $P$, one obtains the original Rudin's Lemma.

\begin{corollary}\label{rudin} \emph{(Rudin's Lemma)} Let $P$ be a poset, $C$ a nonempty lower subset of $P$ and $\mathcal F\in \mathbf{Fin}~P$ a filtered family with $\mathcal F\subseteq\Diamond C$. Then there exists a directed subset $D$ of $C$ such that $\mathcal F\subseteq \Diamond\da D$.
\end{corollary}

For a $T_0$ space $X$ and $\mathcal{K}\subseteq \mathord{\mathsf{K}}(X)$, let $M(\mathcal{K})=\{A\in \mathcal C(X) : K\bigcap A\neq\emptyset \mbox{~for all~} K\in \mathcal{K}\}$ (that is, $\mathcal A\subseteq \Diamond A$) and $m(\mathcal{K})=\{A\in \mathcal C(X) : A \mbox{~is a minimal menber of~} M(\mathcal{K})\}$.

\begin{definition}\label{rudinset} (\cite{Shenchon})
		Let $X$ be a $T_0$ space. A nonempty subset  $A$  of $X$  is said to have the \emph{Rudin property}, if there exists a filtered family $\mathcal K\subseteq \mathord{\mathsf{K}}(X)$ such that $\overline{A}\in m(\mathcal K)$ (that is,  $\overline{A}$ is a minimal closed set that intersects all members of $\mathcal K$). Let $\mathsf{RD}(X)=\{A\in \mathcal C(X) : A\mbox{~has Rudin property}\}$.

The sets in $\mathsf{RD}(X)$ will also be called \emph{Rudin sets}.
\end{definition}
	
The Rudin property is called the \emph{compactly filtered property} in \cite{Shenchon}. In order to emphasize its origin, here we call such property the Rudin property.

\begin{proposition}\label{rudinwf}
	Let $X$ be a $T_0$ space and  $Y$ a well-filtered space. If $f : X\longrightarrow Y$ is continuous and $A\subseteq X$ is a Rudin set, then there exists a unique $y_A\in X$ such that $\overline{f(A)}=\overline{\{y_A\}}$.
\end{proposition}
\begin{proof} Since $A$ has Rudin property, there exists a filtered family $\mathcal K\subseteq \mathord{\mathsf{K}}(X)$ such that $\overline{A}\in m(\mathcal K)$. Let $\mathcal{K}_f=\{\ua f(K\cap \overline{A}) : K\in \mathcal K\}$. Then $\mathcal{F}_f\subseteq \mathord{\mathsf{K}}(Y)$ is filtered. By the proof of Lemma \ref{rudinimage}, $\overline{f(A)}\in m(\mathcal{K}_f)$. By the well-filteredness of $Y$, $\bigcap_{K\in \mathcal{K}}\ua f(K\cap \overline{A})\cap \overline{f(A)}\neq\emptyset$. Select a $y_A\in \bigcap_{K\in \mathcal K} \ua f(K\cap \overline{A})\cap \overline{f(A)}$, then $\overline{\{y_A\}}\subseteq \overline{f(A)}$ and $K\cap \overline{A}\cap f^{-1}(\overline{\{y_A\}})\neq\emptyset$ for all $K\in \mathcal K$. It follows that $\overline{A}=\overline{A}\cap f^{-1}(\overline{\{y_A\}})$ by the minimality of $\overline{A}$, and consequently, $\overline{f(A)}\subseteq \overline{\{y_A\}}$. Therefore, $\overline{f(A)}=\overline{\{y_A\}}$. The uniqueness of $y_A$ follows from the $T_0$ separation of $Y$.
\end{proof}

Motivated by Proposition \ref{rudinwf}, we give the following definition.

\begin{definition}\label{WDspace}
	 A subset $A$ of a space $X$ is called a \emph{well-filtered determined set} if for any continuous mapping $ f:X\longrightarrow Y$
into a well-filtered space $Y$, there exists a unique $y_A\in Y$ such that $\overline{f(A)}=\overline{\{y_A\}}$.
Denote by $\mathsf{WD}(X)$ the set of all closed well-filtered determined subsets of $X$.
\end{definition}

Obviously, a subset $A$ of a space $X$ is well-filtered determined if{}f $\overline{A}$ is well-filtered determined.

\begin{proposition}\label{DRWIsetrelation}
	Let $X$ be a $T_0$ space. Then $S_c(X)\subseteq\mathcal{D}_c(X)\subseteq \mathsf{RD}(X)\subseteq\mathsf{WD}(X)\subseteq\ir_c(X)$.
\end{proposition}
\begin{proof}
	Obviously, $S_c(X)\subseteq\mathcal{D}_c(X)$. Now we prove that the closure of a directed subset $D$ of $X$ is a Rudin set. Let
$\mathcal K_D=\{\ua d : d\in D\}$. Then $\mathcal K_D\subseteq \mathord{\mathsf{K}}(X)$ is filtered and $\overline{D}\in M(\mathcal K_D)$. If $A\in M(\mathcal K_D)$, then $d\in A$ for every $d\in D$, hence $\overline{D}\subseteq A$. So $\overline{D}\in m(\mathcal K_D)$. Therefore $\overline{D}\in \mathsf{RD}(X)$. By Proposition \ref{rudinwf}, we have $\mathsf{RD}(X)\subseteq\mathsf{WD}(X)$. Finally, We show $\mathsf{WD}(X)\subseteq\ir_c(X)$. Let $A\in\mathsf{WD}(X)$.
	Since $\eta_X: X\longrightarrow X^s,\ x\mapsto\da x$, is a continuous mapping into a well-filtered space ($X^s$ is sober), there exists $C\in \ir_c(X)$ such that $\overline{\eta_X(A)}=\overline{\{C\}}$.
	Let $U\in\mathcal O(X)$.
	Note that $$\begin{array}{lll}
	A\cap U\neq\emptyset &\Leftrightarrow& \eta_X(A)\cap  \Diamond U\neq\emptyset\\
	&\Leftrightarrow&\{C\}\cap \Diamond U\neq\emptyset\\
	&\Leftrightarrow& C\in \Diamond U\\
	&\Leftrightarrow& C\cap U\neq\emptyset.
	\end{array}$$
	This implies that $A=C$, so $A\in \ir_c(X)$.
\end{proof}

\begin{lemma}\label{rudinimage}
	Let $X, Y$ be two $T_0$ spaces and $f : X\longrightarrow Y$ a continuous mapping.
\begin{enumerate}[\rm (1)]
	\item If $A\in \mathsf{RD}(X)$, then $\overline{f(A)}\in \mathsf{RD}(Y)$.
	\item If $A\in \mathsf{WD}(X)$, then $\overline{f(A)}\in \mathsf{WD}(Y)$.
\end{enumerate}
\end{lemma}
\begin{proof} (1): It has been proved in \cite{Shenchon}. Here we give a more direct proof. Since $A\in \mathsf{RD}(X)$, there exists a filtered family $\mathcal K\subseteq \mathord{\mathsf{K}}(X)$ such that $A\in m(\mathcal K)$. Let $\mathcal{K}_f=\{\ua f(K\cap A) : K\in \mathcal K\}$. Then $\mathcal{K}_f\subseteq \mathord{\mathsf{K}}(Y)$ is filtered. For each $K\in \mathcal K$, since $K\cap A\neq\emptyset$, we have $\emptyset\neq f(K\cap A)\subseteq \ua f(K\cap A)\cap \overline{f(A)}$. So $\overline{f(A)}\in M(\mathcal{K}_f)$. If $B$ is a closed subset of $\overline{f(A)}$ with $B\in M(\mathcal{K}_f)$, then $B\cap\ua f(K\cap A)\neq\emptyset$ for every $K\in \mathcal K$. So $K\cap A\cap f^{-1}(B)\neq\emptyset$ for all $K\in \mathcal K$. It follows that $A=A\cap f^{-1}(B)$ by the minimality of $A$, and consequently, $\overline{f(A)}\subseteq B$. Therefore, $\overline{f(A)}=B$. Thus $\overline{f(A)}\in \mathsf{RD}(Y)$.

(2): Let $Z$ be a well-filtered space and $g:Y\longrightarrow Z$ a continuous mapping.
Since $g\circ f:X\longrightarrow Z$ is continuous and $A\in \wdd (X)$, there is $z\in Z$ such that $\overline{g(\overline{f(A)})}=\overline{g\circ f(A)}=\overline{\{z\}}$. Thus $\overline{f(A)}\in \wdd (Y)$.
\end{proof}

\begin{lemma}\label{LHCdirect} \emph{(\cite{E_20182})}
	Let $X$ be a locally hypercompact $T_0$ space and $A\in\ir(X)$. Then there exists a directed subset $D\subseteq\da A$ such that $\overline{A}=\overline{D}$.
\end{lemma}

\begin{corollary}\label{LHCdirected}
	For any locally hypercompact  $T_0$ space $X$, $\ir_c (X)=\mathsf{WD}(X)=\mathsf{RD}(X)=\mathcal{D}_c(X)$.
\end{corollary}

\begin{proposition}\label{LCrudin}
	For any locally compact $T_0$ space $X$, $\ir_c(X)=\mathsf{WD}(X)=\mathsf{RD}(X)$.
\end{proposition}
\begin{proof}
	Suppose that $A\in \ir_c(X)$. Let $\mathcal K_A=\{K \in \mk (X) : A\cap\ii \, K\neq\emptyset\}$.
	
	{Claim 1:} $\mathcal K_A\neq\emptyset$.
	
	Let $a\in A$. Since $X$ is locally compact, there exists  $K\in \mk (X)$ such that $a\in\ii \, K$. So $a\in A\cap\ii \, K$ and $K\in \mathcal K_A$, implying $\mathcal{K}_A\neq\emptyset$.
	
	{Claim 2:} $\mathcal K_A$ is filtered.
	
	Let $K_1, K_2\in\mathcal K_A$. Then $A\cap\ii \, K_1\neq\emptyset$ and $A\cap\ii \, K_2\neq\emptyset$. Since $A$ is irreducible, $A\cap\ii \,  K_1\cap\ii \, K_2\neq\emptyset$. Let $x\in A\cap\ii \, K_1\cap\ii \, K_2$. By the local compactness of $X$ again, there exists  $K_3\in \mk (X)$ such that $x\in \ii \, K_3\subseteq K_3\subseteq \ii \, K_1\cap\ii \, K_2$. Thus $K_3\in\mathcal K_A$ and $K_3\subseteq K_1\cap K_2$. So $\mathcal K_A$ is filtered.

    {Claim 3:}  $A\in m(\mathcal K_A)$.	

Clearly, $\mathcal K_A\subseteq \Diamond A$. If $B$ is a proper closed subset of $A$, then there is $a\in A\setminus B$. Since $X$ is locally compact, there is $K_a\in \mk (X)$ such that $a\in \ii \, K_a\subseteq K_a\subseteq X\setminus B$. Then $K_a\in \mathcal K_A$ but $K_a\cap B=\emptyset$. Thus $A$ is a minimal closed set that meets all members of $\mathcal K_A$, and hence $A\in\mathsf{RD}(X)$. By Proposition \ref{DRWIsetrelation}, $\ir_c(X)=\mathsf{WD}(X)=\mathsf{RD}(X)$.
\end{proof}

\section{Some new characterizations of well-filtered spaces and sober spaces}

In \cite{Xi-Lawson-2017}, Xi and Lawson have given a sufficient condition for a $T_0$ space to be well-filtered. We now give some new characterizations of well-filtered and sober spaces based on the results in the above section.

\begin{proposition}\label{wfrudinc}
Let $X$ be a $T_0$ space. Then the following two conditions are equivalent:
	\begin{enumerate}[\rm (1)]
		\item $X$ is well-filtered.
		\item $\wdd (X)=\mathcal S_c(X)$.
        \item $\mathsf{RD}(X)=\mathcal S_c(X)$.
	\end{enumerate}
\end{proposition}
\begin{proof} (1) $\Rightarrow$ (2):  Use the mapping $id_X : X \longrightarrow X$.

(2) $\Rightarrow$ (3): By Proposition \ref{DRWIsetrelation}.

(3) $\Rightarrow$ (1): Suppose that $\mathcal K\subseteq \mathord{\mathsf{K}}(X)$ is filtered, $U\in \mathcal O(X)$, and $\bigcap \mathcal K \subseteq U$. If $K\not\subseteq U$ for each $K\in \mathcal K$, then by Lemma \ref{t Rudin}, the closed set $X\setminus U$ contains an irreducible closed subset $A$ that also meets all members of $\mathcal{K}$ and hence $A\in \mathsf{RD}(X)$. By (2), $A=\overline{\{x\}}$ for some $x\in X$. For each by $K\in \mathcal K$, since $K\cap A=K\cap \overline{\{x\}}\neq\emptyset$, we have $x\in K$. So $x\in \bigcap \mathcal K \subseteq U\subseteq X\setminus A$, a contradiction. Therefore, $K\subseteq U$ for some $K\in \mathcal K$.
\end{proof}

\begin{remark}\label{shenc} The equivalence of (1) and (3) in Proposition \ref{wfrudinc} has been proved in \cite{Shenchon} in a different way.
\end{remark}

By Proposition \ref{d-spaceCharact} and Proposition \ref{wfrudinc}, we get the following result.

\begin{corollary}\label{wf-dspace} Every well-filtered space is a $d$-space.
\end{corollary}

\begin{corollary}\label{wfretract}  A retract of a well-filtered space is well-filtered.
\end{corollary}
\begin{proof} Suppose that $Y$ is a retract of a well-filtered space $X$. Then there are continuous mappings $f : X\longrightarrow Y$ and $g : Y\longrightarrow X$ such that $f\circ g=id_Y$. Let $B\in \mathsf{RD}(Y)$, then by Lemma \ref{rudinimage} and Proposition \ref{wfrudinc}, there exists a unique $x_B\in X$ such that $\overline{g(B)}=\overline{\{x_B\}}$. Therefore, $\overline{B}=\overline{f\circ g(B)}=\overline{f(\overline{g(B)})}=\overline{f(\overline{\{x_B\}})}=\overline{\{f(x_B)\}}$. By Proposition \ref{wfrudinc}, $Y$ is well-filtered.
\end{proof}

By Theorem \ref{LCrudin} and Proposition \ref{wfrudinc}, we get the following well-known result.

\begin{corollary}\label{LCwfissober}\emph{(\cite{redbook}\cite{Kou})}
Every locally compact well-filtered space is sober.
\end{corollary}

\begin{example}\label{examp1}
	Let $X$ be a countably infinite set and $X_{cof}$ the space equipped with the \emph{co-finite topology} (the empty set and the complements of finite subsets of $X$ are open). Then
\begin{enumerate}[\rm (a)]
    \item $\mathcal C(X_{cof})=\{\emptyset, X\}\cup X^{(<\omega)}$, $X_{cof}$ is $T_1$ and hence a strong $d$-space.
    (see Definition \ref{strong $d$-space}).
    \item $\mk (X_{cof})=2^X\setminus \{\emptyset\}$.
    \item $X_{cof}$ is locally compact and first countable.
    \item $\mathcal S_c(X_{cof})=\mathcal D_c(X_{cof})=\{\{x\} : x\in X\}$,
    \item $\ir_c(X_{cof})=\mathrm{WD}(X_{cof})=\mathrm{RD}(X_{cof})=\{X\}\cup \{\{x\} : x\in X\}\neq\mathcal D_c(X_{cof})$. In fact, let $\mathcal{K}_X=\{X\setminus F : F\in X^{(<\omega)}\}$. Then $\mathcal{K}_X\in \mk (X_{cof})$ is filtered and $X\in M(\mathcal{K}_X)$. For any $A\in \mathcal C(X_{cof})$, if $A\neq X$, then $A$ is finite and hence $A\not\in M(\mathcal{K}_X)$ because $A\cap (X\setminus A)=\emptyset$. Thus $X\in \mathrm{RD}(X_{cof})$, but $X\not\in D_c(X_{cof})$.
\item $X_{cof}$ is not well-filtered (by Proposition \ref{wfrudinc}).
\end{enumerate}
\end{example}

\begin{example}\label{examp2}
	Let $L$ be the complete lattice constructed by Isbell \cite{isbell}. Then by \cite[Corollary 3.2]{Xi-Lawson-2017} (or Corollary \ref{Scott is WF2} below), $X=\Sigma~\!\!L$ is a well-filtered space. Note that $\Sigma L$ is not sober. Thus by Proposition \ref{wfrudinc}, $\wdd(X)\neq\ir_c(X)$ and $\kf(X)\neq\ir_c(X)$.
\end{example}

\begin{example}\label{examp3}
	Let $X$ be a uncountably infinite set and $X_{coc}$ the space equipped with \emph{the co-countable topology} (the empty set and the complements of countable subsets of $X$ are open). Then
\begin{enumerate}[\rm (a)]
    \item $\mathcal C(X_{coc})=\{\emptyset, X\}\cup X^{(\leqslant\omega)}$, $X_{coc}$ is $T_1$ and hence a strong $d$-space.
    \item $\mk (X_{coc})=X^{(<\omega)}\setminus \{\emptyset\}$ and $\ii~\!K=\emptyset$ for all $K\in \mk (X_{coc})$.
    \item $X_{coc}$ is not locally compact and not first countable.
    \item $\ir_c(X_{coc})=\{X\}\cup \{\{x\} : x\in X\}$, $\mathrm{WD}(X_{coc})=\mathrm{RD}(X_{coc})=\mathcal D_c(X_{coc})=\mathcal S_c(X_{coc})=\{\{x\} : x\in X\}$. Therefore, $\ir_c(X_{coc})\neq \mathrm{WD}(X_{coc})$.
    \item $X_{coc}$ is well-filtered (by Proposition \ref{wfrudinc}), but it not sober.
\end{enumerate}
\end{example}

\begin{theorem}\label{WFxu-zhao1}
	For any $T_0$ space $X$, the following conditions are equivalent:
 \begin{enumerate}[\rm (1)]
		\item $X$ is well-filtered.
        \item For each $(A, K)\in \wdd(X)\times\mathbf{up} (X)$, $max (A)\neq\emptyset$ and $\downarrow (A\cap K)\in \mathcal C(X)$.
        \item For each $(A, K)\in \kf(X)\times\mathbf{up} (X)$, $max (A)\neq\emptyset$ and $\downarrow (A\cap K)\in \mathcal C(X)$.
        \item For each $(A, K)\in \wdd(X)\times\mk (X)$, $max (A)\neq\emptyset$ and $\downarrow (A\cap K)\in \mathcal C(X)$.
        \item For each $(A, K)\in \kf(X)\times\mk (X)$, $max (A)\neq\emptyset$ and $\downarrow (A\cap K)\in \mathcal C(X)$.
\end{enumerate}
\end{theorem}

\begin{proof}  (1) $\Rightarrow$ (2):  Suppose that $X$ is well-filtered and $(A, K)\in \ir_c(X)\times\mathbf{up} (X)$. Then by Proposition \ref{wfrudinc}, there is  $x\in X$ such that $A=\overline {\{x\}}$, and hence $max (A)=\{x\}\neq\emptyset$. Now we show that $\downarrow\! (A\cap K)=\downarrow\! (\downarrow x\cap K)$ is closed. If $\downarrow\! (\downarrow\! x\cap K)\neq\emptyset$ (i.e., $\downarrow x\cap K\neq\emptyset$), then $x\in K$ since $K$ is saturated (that is, $K$ is an upper set). It follows that $\downarrow\! (\downarrow\! x\cap K)=\downarrow\! x\in \mathcal C(X)$.

(2) $\Rightarrow$ (3) and (4) $\Rightarrow$ (5): By Proposition \ref{DRWIsetrelation}.

(2) $\Rightarrow$ (4) and  (3) $\Rightarrow$ (5): Trivial.

(5) $\Rightarrow$ (1): Suppose that $\mathcal K\subseteq \mathord{\mathsf K}(X)$ is filtered, $U\in \mathcal O(X)$, and $\bigcap \mathcal K \subseteq U$. If $K\not\subseteq U$ for each $K\in \mathcal K$, then by Lemma \ref{t Rudin}, $X\setminus U$ contains a minimal irreducible closed subset $A$ that still meets all members of $\mathcal{K}$, and hence $A\in \kf(X)$. For any $\{K_1, K_2\}\subseteq \mathcal K$, we can find $K_3\in \mathcal K$ with $K_3\subseteq K_1\cap K_2$. It follows that $\downarrow (A\cap K_1)\in \mathcal C(X)$ and $\emptyset \neq A\cap K_3\subseteq \downarrow (A\cap K_1)\cap K_2\neq\emptyset$. By (5) and the minimality of $A$, we have $\downarrow (A\cap K_1)=A$ for all $K_1\in\mathcal K$. Select an $x\in max (A)$. Then for each $K\in \mathcal K$, $x\in  \downarrow (A\cap K)$, and consequently, there is $a_k\in A\cap K$ such that $x\leq a_k$. By the maximality of $x$ we have $x=a_k$. Therefore, $x\in K$ for all $K\in \mathcal K$, and so $x\in \bigcap \mathcal K \subseteq U\subseteq X\setminus A$, a contradiction. Thus $X$ is well-filtered.
\end{proof}

Note that if $X$ is a $d$-space, then by Lemma \ref{d-spacemax}, $max(A)\not=\emptyset$ for every closed set $A$ of $X$. Therefore, by Theorem \ref{WFxu-zhao1}, we obtain the following corollary.

\begin{corollary}\label{WFxu-zhao2}
	For a $d$-space $X$, the following conditions are equivalent:
 \begin{enumerate}[\rm (1)]
		\item $X$ is well-filtered.
        \item For each $(A, K)\in \wdd(X)\times\mathbf{up} (X)$, $\downarrow (A\cap K)\in \mathcal C(X)$.
        \item For each $(A, K)\in \kf(X)\times\mathbf{up} (X)$, $\downarrow (A\cap K)\in \mathcal C(X)$.
        \item For each $(A, K)\in \wdd(X)\times\mk (X)$, $\downarrow (A\cap K)\in \mathcal C(X)$.
        \item For each $(A, K)\in \kf(X)\times\mk (X)$, $\downarrow (A\cap K)\in \mathcal C(X)$.
\end{enumerate}
\end{corollary}

\begin{corollary}\label{xi-lawsonWF}\emph{(\cite{Xi-Lawson-2017})}
	Let $X$ be a $d$-space such that $\downarrow (A\cap K)$ is closed for all
$A\in \mathcal C(X)$ and $K\in \mk (X)$. Then $X$ is well-filtered.
\end{corollary}

\begin{definition}\label{semiclosed} Let $P$ be a poset equipped with a topology. The partial
order is said to be \emph{upper semiclosed} if each $\ua x
$ is closed.
\end{definition}

\begin{definition}\label{upper semicompact}
A topological space $X$ with a partially order is called \emph{upper semicompact}, if $\ua x$ is compact for any $x\in X$, or equivalently, if $\ua x\cap A$ is compact for any $x\in X$ and $A\in \mathcal C(X)$. $X$ is called \emph{weakly upper semicompact} if $\ua x\cap A$ is compact for any $x\in X$ and $A\in \ir_c(X)$.
\end{definition}

\begin{lemma}\label{u-semicloed closed}\emph{(\cite{redbook})}
Let $X$ be a topological space with an upper semiclosed partial order. If $A$
is a compact subset of $X$, then $\da A$ is Scott closed.
\end{lemma}

\begin{lemma}\label{d-spaceScottclosed} Let $X$ be a $T_0$ space such that $\Sigma ~\!\! X$ is a $d$-space. For $A\in \mathcal C(X)$ and $K\in \mk (\Sigma ~\!\!X)$, if $\da (\ua x\cap A)\in \mathcal C(\Sigma~\!\!X)$ for all $x\in X$, then $\da (K\cap A)=\bigcup_{k\in K}\da (\ua k\cap A)\in \mathcal C(\Sigma ~\!\!X)$.
\end{lemma}
\begin{proof} Since $\Sigma ~\!\! X$ is a $d$-space, $X$ is a dcpo. Let $D\in \mathcal D(X)$ such that $D\subseteq \da (K\cap A)$. If $\bigvee D\not\in  \da (K\cap A)$, then for each $k\in K$, $\bigvee D\not\in  \da (\ua k\cap A)$, and hence $\bigcap_{d\in D}\ua d\cap \da (\ua k\cap A)=\emptyset$. For each $k\in K$, since $\Sigma~\!\! X$ is a $d$-space and $\da (\ua k\cap A)\in \sigma (X)$, by Proposition \ref{d-spaceCharact}, there is a $d_k\in D$ such that $\ua d_k\cap \ua k\cap A=\emptyset$, and consequently, $k\in X\setminus \da (\ua d_k \cap A)$ and $\da (\ua d_k \cap A)\in \mathcal C(\Sigma~\!\!X)$. By the compactness of $K$ in $\Sigma ~\!\!X$, there exists a finite subset $\{d_{k_1}, ..., d_{k_n}\}\subseteq D$ such that $K\subseteq \bigcup_{i=1}^{n} (X\setminus \da (\ua d_{k_i} \cap A))$. By the directness of $D$, there is a $d_0$ such that $\ua d_0\subseteq \bigcap_{i=1}^{n}\ua d_{k_i}$. It follows that $K\subseteq X\setminus \da (\ua d_o\cap A)$, which contradicts $d_o\in \da (K\cap A)$, hence $\bigvee D\in  \da (K\cap A)$.
\end{proof}

By Corollary \ref{WFxu-zhao2}, Lemma \ref{u-semicloed closed} and Lemma \ref{d-spaceScottclosed}, we get the following corollaries.

\begin{corollary}\label{Scott is WF1} For a dcpo $P$, if $(P, \lambda (P))$ is weakly upper semicompact, then $(P, \sigma (P))$ is well-filtered.
\end{corollary}

\begin{corollary}\label{Scott is WF2}\emph{(\cite{Xi-Lawson-2017})} For a dcpo $P$, if $(P, \lambda (P))$ is upper semicompact \emph{(}in particular, if  $P$ is bounded complete\emph{)}, then $(P, \sigma (P))$ is well-filtered.
\end{corollary}

In order to reveal finer links between $d$-spaces and well-filtered spaces, we introduce another class of $T_0$ spaces.

\begin{definition}\label{strong $d$-space} A $T_0$ space $X$ is called a strong $d$-space if for any $D\in \mathcal D(X)$, $x\in X$ and $U\in \mathcal O(X)$, $\bigcap_{d\in D}\ua d\cap \ua x\subseteq U$ implies $\ua d\cap \ua x\subseteq U$ for some $d\in D$.
\end{definition}

Clearly, $X$ is a strong $d$-space if{}f for any $D\in \mathcal D(X)$, $x\in X$ and $A\in \mathcal C(X)$, if $\ua d\cap \ua x\cap A\neq\emptyset$ for all $d\in D$, then $\bigcap_{d\in D}(\ua d\cap \ua x)\cap A\neq\emptyset$.

Every $T_1$ space is a strong $d$-space. Also it is easy to verify that every coherent well-filtered space is a strong $d$-space.

\begin{proposition}\label{strong $d$-space by finite upper sets} For a $T_0$ space $X$, the following two conditions are equivalent:
  \begin{enumerate}[\rm (1)]
		\item $X$ is a strong $d$-space.
        \item For any $D\in \mathcal{D}(X)$, $\ua F\in \mathbf{Fin}(X)$ and $U\in \mathcal O(X)$, $\bigcap_{d\in D}\ua d\cap \ua F\subseteq U$ implies $\ua d\cap \ua F\subseteq U$ for some $d\in D$.
\end{enumerate}
\end{proposition}
\begin{proof} (1) $\Rightarrow$ (2): Let $D\in \mathcal D$, $\ua F\in \mathbf{Fin}(X)$ and $U\in \mathcal O(X)$ such that $\bigcap_{d\in D}\ua d\cap \ua F\subseteq U$. Then for each $u\in F$,  $\bigcap_{d\in D}\ua d\cap \ua u\subseteq U$, and hence $\ua d_u\cap \ua u\subseteq U$ for some $d_u\in D$. Since $F$ is finite and $D$ is a direct subset of $X$, there is a $d_0\in D$ such that $\ua d_0\subseteq \bigcap_{u\in F}\ua u$. It follows that $\ua d_0\cap \ua F\subseteq U$.

(2) $\Rightarrow$ (1): Trivial.
\end{proof}

Figure 1 shows certain relations of some spaces lying between $d$-spaces and $T_2$ spaces.

\begin{figure}[ht]
	\centering
	\includegraphics[height=1.5in,width=4.0in]{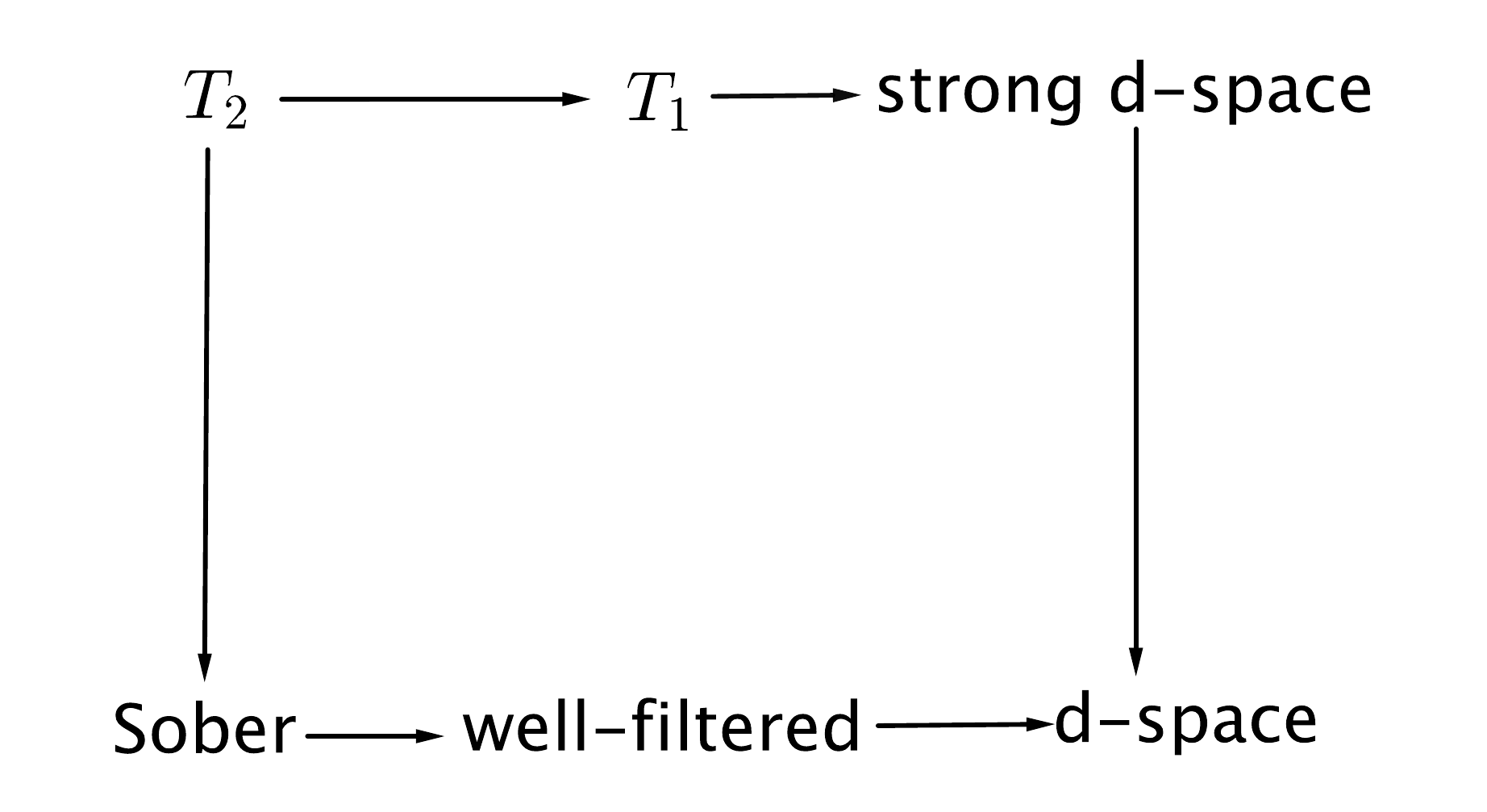}
	\caption{Relations of some spaces lying between $d$-spaces and $T_2$ spaces}
\end{figure}

\begin{definition}\label{strong Scott topology} Let $P$ be a dcpo. A subset $U\subseteq P$ is called \emph{strongly Scott open} if
(i) $U=\mathord{\uparrow}U$, and (ii) for any $D\in \mathcal D(P)$ and $x\in P$, $\bigcap_{d\in D}\ua d\cap \ua x\subseteq U$ (that is, $\ua \bigvee D\cap\ua x\subseteq U$) implies $\ua d\cap \ua x\subseteq U$ for some $d\in D$. Let $\sigma^s (P)$ denote the set of all strongly Scott open subsets of $P$.

Clearly, if $U, V\in \sigma^s (P)$, then $U\cap V\in \sigma^s (P)$. The topology generated by $\sigma^s (P)$ (as a base) is called the \emph{strong Scott topology} on $P$ and
denote it by $\sigma_s(P)$. The space $(P,\sigma_s(P))$ is called the
\emph{strong Scott space} of $P$, and will be denote by $\Sigma_s~\!\! P$.
\end{definition}

For any $x, y$ in a poset $P$, if $\ua x\cap \da y\neq\emptyset$, then $\da (\ua x\cap \da y)=\da y$, and whence for any nonempty finite subset $F$ of $P$, $\da (\ua x\cap \da F)=\emptyset$ or $\da (\ua x\cap \da F)=\da F_x$, where $F_x=\{u\in F :  \ua x\cap \da u\neq\emptyset\}$. Now we check $P\setminus \da F\in \sigma^s (P)$. For any $D\in \mathcal D(P)$ and $x\in X$, if $\bigcap_{d\in D}\ua d\cap \ua x\subseteq P\setminus \da F$, then $\bigcap_{d\in D}\ua d\cap \ua x\cap \da F=\emptyset$, and whence $\bigcap_{d\in D}\ua d\cap \da (\ua x\cap \da F)=\emptyset$, or equivalently, $\bigcap_{d\in D}\ua d\subseteq P\setminus \da (\ua x\cap \da F)\in \upsilon (P)\subseteq \sigma (P)$. Therefore, $\ua d\subseteq P\setminus \da (\ua x\cap \da F)$ for some $d\in D$, and so $\ua d\cap \ua x\subseteq P\setminus \da F$. Thus we have the following relations:
$$\upsilon (P)\subseteq \sigma_s(P)\subseteq\sigma (P).$$
Therefore, if $(P, \upsilon (P))$ is not a strong $d$-space, then any space $(P, \tau)$ with $\upsilon (P)\subseteq \tau \subseteq \alpha (P)$ is not a strong $d$-space. In particular, $\Sigma_s~\!\!P$ and $\Sigma~\!\! P$ are not strong $d$-spaces. For two topologies $\tau_1$ and $\tau_2$ on $P$ with $\upsilon (P)\subseteq \tau_1\subseteq \tau_2$, if $(P, \tau_2)$ is a strong $d$-space, then $(P, \tau_1)$ is also a strong $d$-space.

\begin{remark}\label{sScottsd-space}

\begin{enumerate}[(1)]
\item If a dcpo $P$ is a sup semilattice, then $\sigma_s (P)=\sigma (P)$. In this case, $(P, \sigma (P))$ is a strong $d$-space.

\item For a dcpo $P$, if $\sigma_s (P)=\sigma^s (P)$, then $\Sigma_s~\!\!P$ is a strong $d$-space.
\end{enumerate}
\end{remark}

\begin{proposition}\label{d-space} For a $T_0$ space $X$, consider the following two conditions:
\begin{enumerate}[(1)]
            \item $X$ is a strong $d$-space.
            \item  $X$ is a $d$-space, and $\mathcal O(X)\subseteq \sigma_s(X)$.
\end{enumerate}
Then \emph{(1)} $\Rightarrow$ \emph{(2)}, and the two conditions are equivalent if $X$ \emph{(}with the specialization order\emph{)} is a sup semilattice.
\end{proposition}
\begin{proof} (1) $\Rightarrow$ (2): Let $D\in \mathcal D(X)$ and $U\in \mathcal O(X)$ with $\bigcap_{d\in D}\ua d\subseteq U$. Then for any $c\in D$,
$$\bigcap_{d\in D}\ua d\cap \ua c=\bigcap_{d\in D}\ua d\subseteq U,$$
and hence $\ua d\cap \ua c\subseteq U$ for some $d\in D$. Take one  $e\in D$ with $d\leq e$ and $e\leq c$. Then $\ua e\subseteq \ua d\cap\ua c\subseteq U$. By Proposition \ref{d-spaceCharact}, $X$ is a $d$-space. Now we prove that $\mathcal O(X)\subseteq \sigma_s(X)$. Suppose $U\in \mathcal O(X)$, $x\in X$ and $D\in \mathcal D(X)$ such that $\bigcap_{d\in D}\ua d\cap \ua x\subseteq U$. Since $X$ is a strong $d$-space, $\ua d\cap \ua x\subseteq U$ for some $d\in D$. Thus $U$ is strongly Scott open.

(2) $\Rightarrow$ (1): Suppose that $X$ is a sup semilattice. For any $D\in \mathcal D$, $x\in X$ and $U\in \mathcal O(X)$, if $\bigcap_{d\in D}\ua d\cap \ua x\subseteq U$, then $\bigcap_{d\in D}\ua (d\vee x)\subseteq U$ and $\{d\vee x : d\in D\}\in \mathcal D(X)$. By Proposition \ref{d-spaceCharact}, $\ua d\cap \ua x=\ua d\vee x\subseteq U$ for some $d\in D$. Thus $X$ is a strong $d$-space.
\end{proof}

\begin{definition} A poset $P$ is said to have  \emph{property D} if for any nonempty subset $\{x_i : i\in I\}\subseteq P$ that has a lower bound (that is $\bigcap_{i\in I}\da x_i\neq\emptyset$), $\bigcap_{i\in I}\da x_i\in \mathrm{Id} (P)$.
\end{definition}

 Clearly, every bounded complete poset has property $D$. For a dcpo $P$, $P$ satisfies property $D$ if{}f every nonempty subset $\{x_i : i\in I\}\subseteq P$ that has a lower bound has the greatest lower bound (that is, $\bigcap_{i\in I}\da x_i$ is a principal ideal of $P$).

 \begin{lemma}\label{uppertop-sd-space} For a poset $P$, the following two conditions are equivalent:
 \begin{enumerate}[\rm (1)]
 \item $(P, \upsilon (P))$ is a strong $d$-space.
 \item $P$ is a dcpo, and for any $\{F_i : i\in I\}\subseteq P^{(<\omega)}$ and $x\in P$, $\da (\ua x\cap \bigcap_{i\in I} \da F_i)\in \mathcal C (\Sigma~\!\! P)$.
 \end{enumerate}
 \end{lemma}
 \begin{proof} (1) $\Rightarrow$ (2): Suppose that $(P, \upsilon (P))$ is a strong $d$-space. Then $(P, \upsilon (P))$ is a $d$-space, and hence $P$ is a dcpo. For $\{F_i : i\in I\}\subseteq P^{(<\omega)}$ and $x\in P$, we show that $\da \bigcap_{i\in I} (\ua x\cap \da F_i)\in \mathcal C (\Sigma~\!\! P)$. For any  $D\in \mathcal D(P)$ with $D\subseteq \da (\ua x\cap \bigcap_{i\in I} \da F_i)$, if $\bigvee D\not\in \da (\ua x\cap \bigcap_{i\in I} \da F_i)$, then $\ua \bigvee D\cap \ua x\cap \bigcap_{i\in I} \da F_i=\bigcap_{d\in D}\ua d\cap \ua x\cap \bigcap_{i\in I} \da F_i=\emptyset$, and whence $\bigcap_{d\in D}\ua d\cap \ua x\subseteq P\setminus \bigcap_{i\in I} \da F_i\in \upsilon (P)$. Since $(P, \upsilon (P))$ is a strong $d$-space, there is $d\in D$ such that $\ua d\cap
 \ua x\subseteq P\setminus \bigcap_{i\in I} \da F_i$, which is a contradiction with $d\in \da (\ua x\cap \bigcap_{i\in I} \da F_i)$. Therefore, $\bigvee D\in \da(\ua x\cap \bigcap_{i\in I} \da F_i)$. Thus $\da (\ua x\cap \bigcap_{i\in I} \da F_i)\in \mathcal C(\Sigma~\!\! P)$.

 (2) $\Rightarrow$ (1): For any $D\in \mathcal D(P)$, $x\in P$ and $U\in \upsilon (P)$ such that $\bigcap_{d\in D}\ua d\cap \ua x\subseteq U$, if $U=P$, then $\ua d\cap \ua x\subseteq U$ for all $d\in D$; if $U$ is a proper $\upsilon$-open subset of $P$, then there is a family $\{F_i : i\in I\}\subseteq P^{(<\omega)}$ such that $U=P\setminus \bigcap_{i\in I}\da F_i$. Therefore, $\ua \bigvee D=\bigcap_{d\in D}\ua d\subseteq P\setminus \da (\ua x\cap\bigcap_{i\in I}\da F_i)\in \sigma (P)$ by condition (2). It follows that $\ua d\subseteq P\setminus \da (\ua x\cap\bigcap_{i\in I}\da F_i)$ for some $d\in D$, and whence $\ua d\cap \ua x\subseteq P\setminus \bigcap_{i\in I}\da F_i$, proving that $(P, \upsilon (P))$ is a strong $d$-space.
\end{proof}

Similarly, we have the following result.

 \begin{lemma}\label{Sctttop-sd-space} For a poset $P$, the following two conditions are equivalent:
 \begin{enumerate}[\rm (1)]
 \item $\Sigma~\!\!P$ is a strong $d$-space.
 \item $P$ is a dcpo, and for any $A\in \mathcal C(\Sigma~\!\!P)$ and $x\in P$, $\da (\ua x\cap A)\in \mathcal C(\Sigma~\!\!P)$.
 \end{enumerate}
 \end{lemma}

 \begin{corollary}\label{uppertopologSd} For a dcpo $P$ satisfying property $D$ \emph{(}in particular, $P$ is a complete semilattice\emph{)}, $(P, \upsilon (P))$ is a strong $d$-space.
 \end{corollary}

\begin{proof} For any $\{F_i : i\in I\}\subseteq P^{(<\omega)}$ and $x\in P$, we show that $\da (\ua x\cap \bigcap_{i\in I}\da F_i)=\bigcap_{i\in I}\da (\ua x\cap\da F_i)$. Obviously, $\da (\ua x\cap \bigcap_{i\in I}\da F_i)\subseteq\bigcap_{i\in I}\da (\ua x\cap\da F_i)$. Conversely, if $y\in \bigcap_{i\in I}\da (\ua x\cap\da F_i)$, then for each $i\in I$, there exists $u_i\in \ua x\cap \da F_i$ with $y\leq u_i$, and hence there is $t_i\in F_i$ such that $u_i\leq t_i$. It follows that $x, y\in \bigcap_{i\in I}\da t_i$. Since $P$ satisfies property $D$, there is a $z\in \bigcap_{i\in I}\da t_i\subseteq \bigcap_{i\in I}\da F_i$ such that $y\leq z$ and $x\leq z$, and whence $z\in \ua x\cap \bigcap_{i\in I}\da F_i$ and $y\in \da (\ua x\cap \bigcap_{i\in I}\da F_i)$. Thus $\da (\ua x\cap \bigcap_{i\in I}\da F_i)=\bigcap_{i\in I}\da (\ua x\cap\da F_i)$. For any $s, t\in P$, if $\ua s\cap \da t\neq\emptyset$, then $t\in \ua s\cap \da t$, and hence $\da (\ua s\cap \da t)=\da t$. Therefore, for each $i\in I$, $\da (\ua x\cap\da F_i)=\bigcup_{t\in F_i}\da (\ua x\cap\da t)=\da F_i^x$, where $F_i^x=\{t\in F_i :  \ua x\cap \da t\neq\emptyset\}$. It follows that $\da (\ua x\cap \bigcap_{i\in I}\da F_i)=\bigcap_{i\in I}\da (\ua x\cap\da F_i)\in \mathcal C((P, \upsilon (P)))\subseteq \mathcal C(\Sigma~\!\! P)$. By lemma \ref{uppertop-sd-space}, $(P, \upsilon (P))$ is a strong $d$-space.
\end{proof}

\begin{remark} For a dcpo $P$, consider the following three conditions:
\begin{enumerate}[\rm (1)]
\item $P$ has property $D$.
\item For any $\{F_i : i\in I\}\subseteq P^{(<\omega)}$ and $x\in P$, $\da (\bigcap_{i\in I}(\ua x\cap \da F_i))=\bigcap_{i\in I}\da (\ua x\cap\da F_i)$ (note that $\bigcap_{i\in I}\da (\ua x\cap\da F_i)$ is always $\upsilon$-closed).
    \item $(P, \upsilon (P))$ is a $d$-space.
\end{enumerate}
Then by Lemma \ref{uppertop-sd-space} and the proof of Corollary \ref{uppertopologSd}, we have (1) $\Rightarrow$ (2) $\Rightarrow$ (3).
\end{remark}

\begin{proposition}\label{d-spacestrong d-space} If $X$ is a $d$-space and $\da (\ua x\cap A)\in \mathcal C(X)$ for all $x\in X$ and $A\in \mathcal C(X)$, then $X$ is a strong $d$-space.
\end{proposition}
\begin{proof} Suppose that $D\in \mathcal D$, $x\in X$ and $U\in \mathcal O(X)$ such that $\bigcap_{d\in D}\ua d\cap \ua x\subseteq U$. Let $A=X\setminus U$. Then $A\in \mathcal C(X)$. If for any $d\in D$, $\ua d\cap \ua x\not\subseteq U$, then $\ua d\cap \da (\ua x\cap A)\neq\emptyset$. Since $X$ is a $d$-space and $\da (\ua x\cap A)\in \mathcal C(X)$, by Proposition \ref{d-spaceCharact}, we have $\bigcap_{d\in D}\ua d\cap \da (\ua x\cap A)\neq\emptyset$, and hence $\bigcap_{d\in D}\ua d\cap \ua x\cap A\neq\emptyset$, a contradiction. Thus $\ua d\cap \ua x\subseteq U$ for some $d\in D$.
\end{proof}

By Theorem \ref{WFxu-zhao1}, Corollary \ref{WFxu-zhao2}, Lemma \ref{u-semicloed closed} and Proposition \ref{d-spacestrong d-space}, we get the following corollaries.

\begin{corollary}\label{WFsd-pace1}
	For a $T_0$ space $X$, if for each $B\in \wdd(X)$ and $(A, K)\in \mathcal C(X)\times\mk (X)$, $max (B)\neq\emptyset$ and $\downarrow (A\cap K)\in \mathcal C(X)$, then $X$ is a well-filtered strong $d$-space.
\end{corollary}

\begin{corollary}\label{WFsd-pace2}
	For a $T_0$ space $X$, if for each $B\in \mathrm{RD}(X)$ and $(A, K)\in \mathcal C(X)\times\mk (X)$, $max (B)\neq\emptyset$ and $\downarrow (A\cap K)\in \mathcal C(X)$, then $X$ is a well-filtered strong $d$-space.
\end{corollary}

\begin{corollary}\label{WFsd-pace3}
	Let $X$ be a $d$-space such that $\downarrow (A\cap K)$ is closed for all
$A\in \mathcal C(X)$ and $K\in \mk (X)$. Then $X$ is a well-filtered strong $d$-space.
\end{corollary}

\begin{corollary}\label{Scottstrong d-space} For a dcpo $P$, if $(P, \lambda (P))$ is upper semicompact, then $(P, \sigma (P))$ is a strong $d$-space.
\end{corollary}

The following example shows that for a dcpo $P$, $(P, \upsilon (P))$ and $(P, \sigma (P))$ need not be strong $d$-spaces, although they are always  $d$-spaces.

\begin{example}\label{d-space not strong d-space} (Johnstone space) Let $\mathbb{J}=\mathbb{N}\times (\mathbb{N}\cup \{\omega\})$ with ordering defined by $(j, k)\leq (m, n)$ if{}f $j = m$ and $k \leq n$, or $n =\omega$ and $k\leq m$ (see \cite{johnstone-81}). Then $\mathbb{J}$ is a dcpo, and hence the Johnstone space $\Sigma~\!\mathbb{J}$ is a $d$-space. However, $\Sigma~\!\mathbb{J}$ is not well-filtered (see \cite[Exercise 8.3.9]{Jean-2013}), and hence non-sober. Clearly, $\bigcap_{n\in \mathbb{N}}(\ua (1, n)\cap \ua (2,1))=\emptyset$, but $\ua (1, n)\cap \ua (2,1)=\{(m, \omega): n\leq m\}\neq\emptyset$ for all $n$. Hence $(\mathbb{J}, \upsilon (\mathbb{J}))$ and $\Sigma~\!\mathbb{J}$ are not strong $d$-spaces.
\end{example}

The following example shows that even for a  continuous dcpo $P$, if $\Sigma~\!\!P$ is not coherent,  $\Sigma~\!\!P$ and $\Sigma_s~\!\!P$ may not be strong $d$-spaces.

\begin{example}\label{sober not strong d-space}  Let $C=\{a_1, a_2, ..., a_n, ...\}\cup \{\omega_0\}$ and $P=C\cup\{b\}\cup\{\omega_1, ..., \omega_n, ...\}$ with the order generated  by
\begin{enumerate}[\rm (a)]
            \item $a_1<a_2<...<a_n<a_{n+1}<...$;
            \item  $a_n<\omega_0$ for all $n\in \mathbb{N}$;
            \item $b<\omega_n$ and $a_m<\omega_n$ for all $n, m\in \mathbb{N}$ with $m\leq n$.
\end{enumerate}

\begin{figure}[ht]
	\centering
	\includegraphics[height=2.2in,width=4.0in]{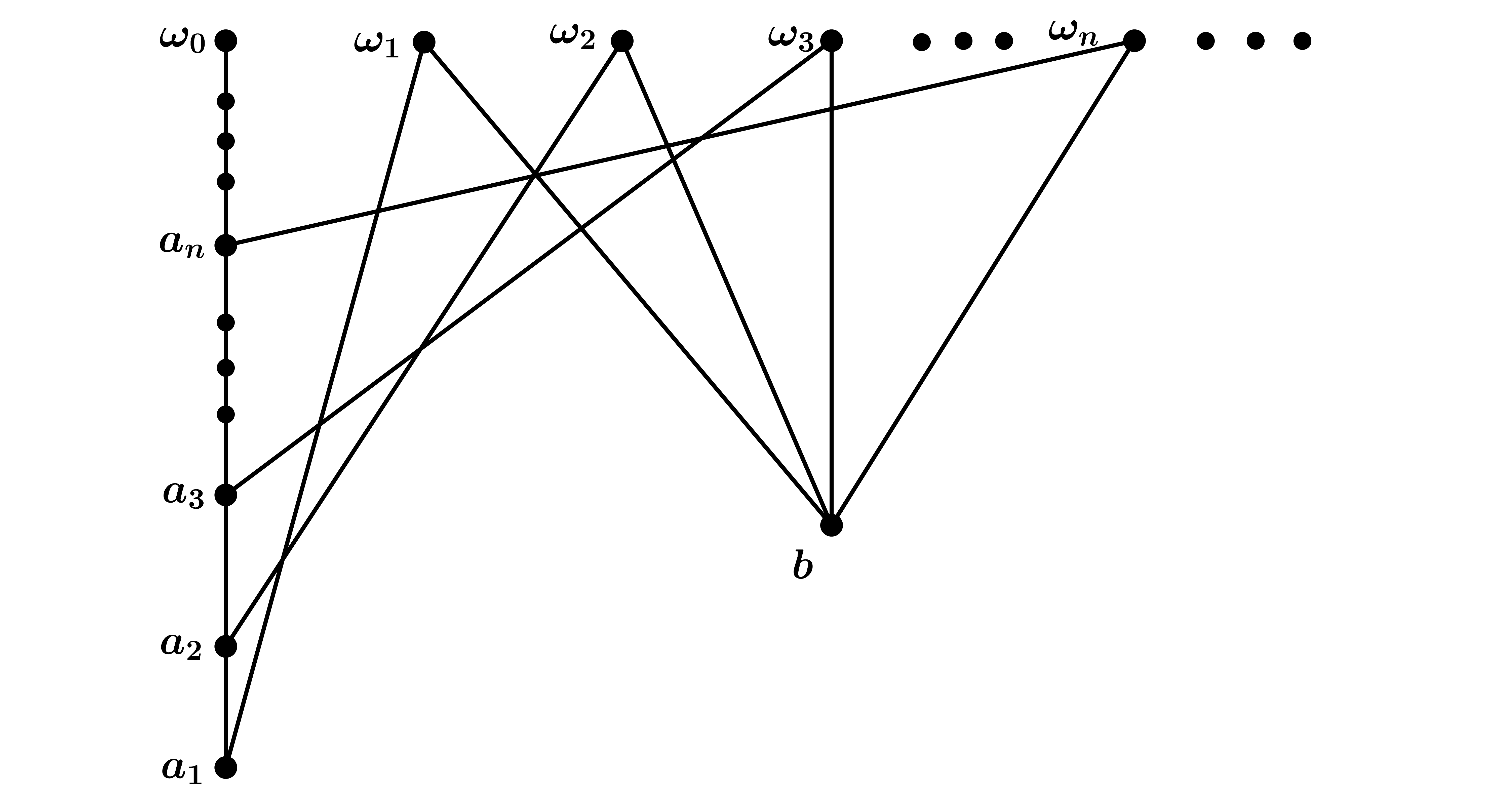}
	\caption{The Poset $P$}
\end{figure}

Then $P$ is a dcpo and $D\in \mathcal D(P)$ if{}f $D\subseteq C$ or $D$ has a largest element (that is $\da D$ is a principal ideal of $P$), and hence $x\ll x$ for all $x\in P\setminus \{\omega_0\}$. Therefore, $P$ is a continuous dcpo and $\Sigma~\!\!P$ is sober. $\ua a_1, \ua b\in \mk (\Sigma~)$, but $\ua a_1\cap \ua b=\{\omega_1, \omega_2, ..., \omega_n, ...\}$ is not Scott compact (note $\{\omega_n\}\in \sigma (P)$ for all $n\in \mathbb{N}$). Thus  $\Sigma~\!\!P$ is not coherent. For any $m\in \mathbb{N}$, $\bigcap_{n\in \mathbb{N}}\ua a_n\cap \ua b=\emptyset \subseteq \{\omega_m\}\in \sigma (P)$, but $\ua a_n\cap \ua b=\{\omega_n, \omega_{n+1}, ...\}\not\subseteq \{\omega_m\}$ for any $n\in \mathbb{N}$. Thus $\Sigma~\!\!P$ is not a strong $d$-space. Since $\bigcap_{n\in \mathbb{N}}\ua a_n\cap \ua b=\emptyset$, but $\ua a_n\cap \ua b=\{\omega_n, \omega_{n+1}, ...\}\neq\emptyset$, $(P, \upsilon (P))$ is not a strong $d$-space, and hence $\Sigma_s~\!\! P$ is not a strong $d$-space because $\upsilon (P)\subseteq\sigma_s (P)$.
\end{example}

\begin{remark} Let $P$ be the continuous domain in Example \ref{sober not strong d-space}. Then $\da (\ua b\cap P)=P\setminus \{\omega_0\}$ is not Scott closed since $\{\omega_0\}\not\in \sigma (P)$ (or equivalently, $\omega_0$ is not a compact element of $P$). So as a sufficient condition for a $d$-space to be well-filtered, the condition that $\downarrow (A\cap K)$ is closed for all
$A\in \mathcal C(X)$ and $K\in \mk (X)$ seems a little too strong (see Corollary \ref{WFxu-zhao2}, Corollary \ref{xi-lawsonWF} and Corollary \ref{WFsd-pace3}).
\end{remark}

\begin{lemma}\label{Smyth-sober}\emph{(\cite{Klause-Heckmann})} For a $T_0$ space $X$, the following conditions are equivalent:
\begin{enumerate}[\rm (1)]
            \item $X$ is sober.
            \item  For any $\mathcal A\subseteq \ir (P_S(X))$ and $U\in \mathcal O(X)$, $\bigcap \mathcal A\subseteq U$ implies $K\subseteq U$ for some $K\in U$.
            \item $P_S(X)$ is sober.
\end{enumerate}
\end{lemma}

\begin{remark}\label{soberH-Keimel} By Remark \ref{meet-in-Smyth} and Lemma \ref{Smyth-sober}, we have that  $X$ is sober if{}f  for any $\mathcal A\subseteq \ir_c (P_S(X))$ and $U\in \mathcal O(X)$, $\bigcap \mathcal A\subseteq U$ implies $K\subseteq U$ for some $K\in\mathcal A$.
\end{remark}

For a $T_0$ space and $A\subseteq X$, define $\Psi_{\mathbf{up}(X)}(A)=\{K\in \mathbf{up}(X) : K\cap A\neq \emptyset\}$ and $\Psi_{\mk (X)}(A)=\{K\in \mk (X) : K\cap A\neq \emptyset\}$.

The following theorem provides a new characterization of sober spaces similar to that for well-filtered spaces given in Theorem \ref{WFxu-zhao1}.

\begin{theorem}\label{soberxu-zhao1}
The following conditions are equivalent for a $T_0$ space $X$:
 \begin{enumerate}[\rm (1)]
		\item $X$ is sober.
        \item For any $(A, K)\in \ir_c(X)\times\mathbf{up}(X)$, $\Psi_{\mathbf{up}(X)}(A)\in \mathrm{Filt}(\mathbf{up}(X))$, $max (A)\neq\emptyset$ and $\downarrow (A\cap K)\in \mathcal C(X)$.
        \item For any $(A, K)\in \ir_c(X)\times\mk (X)$, $\Psi_{\mk (X)}(A)\in \mathrm{Filt}(\mk (X))$, $max (A)\neq\emptyset$ and $\downarrow (A\cap K)\in \mathcal C(X)$.
\end{enumerate}
\end{theorem}

\begin{proof}  We only prove the equivalence of (1) and (3). The proof of the equivalence of (1) and (3) is similar.

(1) $\Rightarrow$ (3):  Suppose that $X$ is sober and $(A, K)\in \ir_c(X)\times\mk (X)$. Then there is an $x\in X$ such that $A=\overline {\{x\}}$, and hence $max (A)=\{x\}\neq\emptyset$. Clearly, $\Psi_{\mk (X)}(A)=\{K\in \mk (X) : K\cap \da x\neq \emptyset\}=\ua_{\mk (X)}\ua x$ is a principal filter of $\mk (X)$. Now we show that $\downarrow (A\cap K)$ is a closed subset of $X$. Obviously, if $\downarrow x\cap K=\emptyset$, then $\downarrow (A\cap K)=\downarrow (\downarrow x\cap K)=\emptyset$; if $\downarrow x\cap K\neq\emptyset$ (in this case $x\in K$ since $K$ is an upper set), then $\downarrow (A\cap K)=\downarrow (\downarrow x\cap K)=\downarrow x$. Thus $\downarrow (A\cap K)=\downarrow x\in \mathcal C(X)$.

(3) $\Rightarrow$ (1): Suppose $\mathcal A\subseteq \ir (P_S(X))$ and $U\in \mathcal O(X)$ such that $\bigcap \mathcal A\subseteq U$. If $K\not\subseteq U$ for all $K\in U$, then by Lemma \ref{t Rudin}, $X\setminus U$ contains a minimal irreducible closed subset $A$ that still meets all members of $\mathcal{A}$. For any $\{K_1, K_2\}\subseteq \mathcal A$, since $\Psi_{\mk (X)}(A)\in \mathrm{Filt}(\mk (X))$, there is a $K_3\in\Psi_{\mk (X)}(A)$ with $K_3\subseteq K_1\cap K_2$. It follows that $\emptyset \neq A\cap K_3\subseteq \downarrow (A\cap K_1)\cap K_2\neq\emptyset$. By (3) and the minimality of $A$, we have $\downarrow (A\cap K_1)=A$ for all $K_1\in\mathcal A$. Select an $x\in max (A)$. Then for each $K\in \mathcal A$, $x\in  \downarrow (A\cap K)$, and consequently, there is $a_k\in A\cap K$ such that $x\leq a_k$. By the maximality of $x$, we have $x=a_k$. Therefore, $x\in K$ for all $K\in \mathcal A$, and whence $x\in \bigcap \mathcal A \subseteq U\subseteq X\setminus A$, a contradiction. By Lemma \ref{Smyth-sober}, $X$ is sober.
\end{proof}

By Lemma \ref{d-spacemax} and Theorem \ref{soberxu-zhao1}, we get the following corollary.

\begin{corollary}\label{soberxu-zhao2}
	For a $d$-space $X$, the following conditions are equivalent:
 \begin{enumerate}[\rm (1)]
		\item $X$ is sober.
        \item For each $(A, K)\in \ir_c(X)\times\mathbf{up}(X)$, $\Psi_{\mathbf{up}(X)}(A)\in \mathrm{Filt}(\mathbf{up}(X))$ and $\downarrow (A\cap K)\in \mathcal C(X)$.
        \item For each $(A, K)\in \ir_c(X)\times\mk (X)$, $\Psi_{\mk (X)}(A)\in \mathrm{Filt}(\mk (X))$ and $\downarrow (A\cap K)\in \mathcal C(X)$.
\end{enumerate}
\end{corollary}

\section{Conclusion}
In this paper, based on topological Rudin's Lemma, we instigated two new classes of subsets lying between the classes of all closures of directed subsets and that of irreducible closed sets. Using such subsets, we obtained some new characterizations of sober spaces and well-filtered spaces, which improve and generalize the related results of Xi and Lawson \cite{Xi-Lawson-2017}.

Our study also leads to the definition of a new class of spaces - strong $d$-spaces and a new topology - the strong Scott topology which may deserve further investigation.

We now close our paper with some problems.
\begin{problem} Does $\mathsf{RD}(X)=\mathsf{WD}(X)$ hold for every $T_0$ space $X$?
\end{problem}

\begin{problem} Does $\ir_c(X)=\mathsf{WD}(X)$ hold for every core compact $T_0$  space $X$?
\end{problem}

\begin{problem} Are the following conditions equivalent for any  $T_0$ space $X$?
\begin{enumerate}[\rm (a)]
\item $X$ is sober.
\item  For each $(A, K)\in \ir_c(X)\times\mathbf{up}(X)$, $max (A)\neq\emptyset$ and $\downarrow (A\cap K)\in \mathcal C(X)$.
\item For each $(A, K)\in \ir_c(X)\times\mk (X)$, $max (A)\neq\emptyset$ and $\downarrow (A\cap K)\in \mathcal C(X)$.
\end{enumerate}
\end{problem}

\begin{problem} Are the following conditions equivalent for any $d$-space $X$?
\begin{enumerate}[\rm (a)]
\item  $X$ is sober.
\item  For each $(A, K)\in \ir_c(X)\times\mathbf{up}(X)$, $\downarrow (A\cap K)\in \mathcal C(X)$.
\item For each $(A, K)\in \ir_c(X)\times\mk (X)$, $\downarrow (A\cap K)\in \mathcal C(X)$.
\end{enumerate}
\end{problem}

\end{document}